\documentclass[11pt, reqno]{amsart}

\title{Tauberian rates and the Baire property of ideals}


\usepackage[T1]{fontenc}
\usepackage{amsmath}
\usepackage{amssymb}
\usepackage{amsthm}
\usepackage[left=3.5cm, right=3.5cm, paperheight=11.8in]{geometry}
\usepackage{hyperref}
\usepackage{fancyhdr}
\usepackage{enumitem}
\usepackage{comment}
\usepackage{tikz}
\usepackage{bm}
\usepackage{nicefrac}
\usepackage{mathrsfs}
\usepackage{graphicx}
\usepackage[utf8]{inputenc}
\usepackage{cancel}
\usepackage{mathtools}

\usetikzlibrary{matrix,arrows,positioning,automata,shapes,shadows,calc,through,intersections}

\AtBeginDocument{%
   \def\MR#1{}
}

\newtheorem{thm}{Theorem}[section]
\newtheorem{cor}[thm]{Corollary}
\newtheorem{lem}[thm]{Lemma}
\newtheorem{prop}[thm]{Proposition}

\theoremstyle{definition} 
\newtheorem{defi}[thm]{Definition}
\let\olddefi\defi
\renewcommand{\defi}{\olddefi\normalfont}
\newtheorem{example}[thm]{Example}
\let\oldexample\example
\renewcommand{\example}{\oldexample\normalfont}
\newtheorem{rmk}[thm]{Remark}
\let\oldrmk\rmk
\renewcommand{\rmk}{\oldrmk\normalfont}

\hypersetup{
    pdftitle={Tauberian rates and the Baire property of ideals},
    pdfauthor={Paolo Leonetti},
    pdfmenubar=false,
    pdffitwindow=true,
    pdfstartview=FitH,
    colorlinks=true,
    linkcolor=blue,
    citecolor=green,
    urlcolor=cyan
}

\providecommand{\MR}[1]{}

\providecommand{\MR}{\relax\ifhmode\unskip\space\fi MR }

\providecommand{\href}[2]{#2}

\author[P.~Leonetti]{Paolo Leonetti}
\address{
Universit\'{a} degli Studi dell’Insubria\\ via Monte Generoso 71 \\ Varese 21100\\ Italy}
\email{leonetti.paolo@gmail.com}
                                 
\begin{document}

\maketitle
\thispagestyle{empty}

\begin{abstract}
We study Tauberian conditions on the successive differences of real sequences which ensure that: (i) ordinary convergence coincides with ideal convergence; or (ii) every ordinary accumulation point is an ideal cluster point; or (iii) every ordinary accumulation point is an ideal limit point. 

First, we prove that an ideal $\mathcal{I}$ has the Baire property if and only if there exists a positive sequence $(r_n)$ tending to zero, with $\sum_n r_n=\infty$, such that 
$$
|x_{n+1}-x_n|=O(r_n)
$$
is a Tauberian condition for any of the three previous statements.  
Second, for analytic $P$-ideals, we characterize the validity of the Tauberian rate $O(1/n)$ in terms of properties of lower semicontinuous submeasures which represent it. 

Lastly, in the case of statistical limit points, we prove that $o(1/n)$ is a Tauberian condition for the statement (iii), whereas $O(1/n)$ need not suffice: in this sense, $o(1/n)$ is an optimal Tauberian rate for (iii). On the other hand, for statistical cluster points in (ii), $O(1/n)$ is an optimal one among nonincreasing rates. 
For every positive null sequence $(r_n)$ with divergent sum, we also construct a summable ideal $\mathcal{I}$ for which $o(r_n)$ is an optimal Tauberian rate for both statements (ii) and (iii).
\end{abstract}

\section{Introduction}\label{sec:intro}

The aim of this work is to provide 
Tauberian conditions on the growth rate of a real sequence $(x_{n})_{n \in \omega}$ which ensure that each of its accumulation points has to be of a rather special type. As a motivating example, we show that if $|x_{n+1}-x_n|=o(1/n)$ as $n\to \infty$ and $\eta$ is an accumulation point of $x$, then there exists a subsequence $(x_n: n \in S)$ convergent to $\eta$ such that, in addition, the index set $S$ does \emph{not} have asymptotic density zero (so that, intuitively, the subsequence has \textquotedblleft sufficiently many elements\textquotedblright). Moreover, the same claim fails if we replace the growth rate $o(1/n)$ with $O(1/n)$, hence this is the best one can hope for; see Theorem \ref{thm:conjectureweak} and Example \ref{ex:nonequivalence} below. 

To be more precise about our main results, we need to recall the following notions. 
Let $\mathcal{I}$ be an \emph{ideal} on the nonnegative integers $\omega$, that is, a family of subsets of $\omega$ which is stable under subsets and finite unions. Unless otherwise stated, we assume that $\{n\} \in \mathcal{I}$ for all $n \in \omega$, and that $\omega\notin \mathcal{I}$. 
The family of $\mathcal{I}$-positive sets is denoted by $\mathcal{I}^+:=\mathcal{P}(\omega)\setminus \mathcal{I}$, while we write $\mathcal{I}^\star:=\{\omega\setminus S: S \in \mathcal{I}\}$ for its dual filter. 
Regarding ideals of subsets of the Cantor space $\{0,1\}^\omega$, it is possible to speak about the topological complexity, e.g., $F_\sigma$ ideals, analytic ideals, etc. Examples of ideals include the family of finite sets $\mathrm{Fin}:=\{S\subseteq \omega: |S|<\infty\}$ (which is an $F_\sigma$ ideal) or the family 
\begin{equation}\label{eq:definitionZ}
\mathcal{Z}:=\left\{S\subseteq \omega: \lim_{n\to \infty} \frac{|S \cap [0,n]|}{n+1}=0\right\}
\end{equation}
of asymptotic density zero sets (which is an $F_{\sigma\delta}$ ideal but not $F_\sigma$). We refer to \cite{MR1711328} for an excellent monograph on ideals (or, equivalently, filters). 

Given a real sequence $x=(x_n)_{n \in \omega}$, we say that $x$ is $\mathcal{I}$\emph{-convergent to} $\eta \in \mathbb{R}$, shortened as $\mathcal{I}\text{-}\lim x=\eta$, if $\{n \in \omega: |x_n-\eta|>\varepsilon\} \in \mathcal{I}$ for all $\varepsilon>0$. Of course, if $\mathcal{I}=\mathrm{Fin}$, this corresponds to ordinary convergence. 
%
In addition, we say that $\eta\in \mathbb{R}$ is an $\mathcal{I}$\emph{-cluster point of} $x$ if $\{n \in \omega: |x_n-\eta|\le \varepsilon\} \in \mathcal{I}^+$ for all $\varepsilon>0$. The set of $\mathcal{I}$-cluster points of $x$ is denoted by $\Gamma_x(\mathcal{I})$. If $\mathcal{I}=\mathrm{Fin}$ then we write $\mathrm{L}_x:=\Gamma_x(\mathcal{I})$. 
%
Lastly, we say that $\eta\in \mathbb{R}$ is an $\mathcal{I}$\emph{-limit point of} $x$ if there exists a subsequence $(x_n: n \in S)$ convergent to $\eta$ such that $S \in \mathcal{I}^+$. The set of $\mathcal{I}$-limit points of $x$ is denoted by $\Lambda_x(\mathcal{I})$. In the literature, if $\mathcal{I}=\mathcal{Z}$, then $\mathcal{Z}$-convergence, $\mathcal{Z}$-cluster points, and $\mathcal{Z}$-limit points are usually referred to as \textquotedblleft statistical convergence,\textquotedblright\, \textquotedblleft statistical cluster points,\textquotedblright\, and \textquotedblleft statistical limit points,\textquotedblright\, respectively; see e.g. \cite{MR4600193, MR1372186, MR1181163, MR1416085}. 



For each real sequence $x$, we define the sequence of differences $\Delta x:=(\Delta x_n)_{n \in \omega}$ by 
$$
\forall n \in \omega, \quad 
\Delta x_n:=x_{n+1}-x_n. 
$$
It is clear that every $\mathcal{I}$-limit point is an $\mathcal{I}$-cluster point and that every $\mathcal{I}$-cluster point is an accumulation point, namely, $\Lambda_x(\mathcal{I})\subseteq \Gamma_x(\mathcal{I})\subseteq \mathrm{L}_x$ for every real sequence $x$. 
Here, we are going to provide Tauberian conditions on the growth rate of $\Delta x$ which ensure that the converse inclusions hold, so that for instance 
$\Gamma_x(\mathcal{I})=\mathrm{L}_x$ or $\Lambda_x(\mathcal{I})=\mathrm{L}_x$. Of course, if $\mathcal{I}=\mathrm{Fin}$, then $\Lambda_x(\mathrm{Fin})=\Gamma_x(\mathrm{Fin})=\mathrm{L}_x$, so that our results will be meaningful whenever $\mathcal{I}$ contains at least one infinite set. A more precise observation follows in the next remark (hereafter, Landau's symbols $O(\cdot)$ and $o(\cdot)$ are always meant implicitly as $n\to \infty$): 
\begin{rmk}\label{rmk:IFin}
For an ideal $\mathcal{I}$ on $\omega$, we have $\mathcal{I}\neq \mathrm{Fin}$ if and only if there exists a real sequence $x$ such that $|\Delta x_n|=O(1)$ and $\Lambda_x(\mathcal{I})=\Gamma_x(\mathcal{I}) \neq \mathrm{L}_x$. The \textsc{If} part is obvious. For the converse, it is enough to consider the sequence $\bm{1}_A$ for some infinite $A \in \mathcal{I}$. 
\end{rmk}

It is worth noting that real sequences $x$ for which $|\Delta x_n|=o(1)$ are usually called \textquotedblleft slow sequences\textquotedblright, see e.g. \cite{MR3785882}. 
Our starting point is the following result due to Fridy, see \cite[Theorem 3]{MR816582}: 
\begin{thm}
\label{thm:fridythm3}
Let $x$ be a real sequence such that $|\Delta x_n|=O(1/n)$. 
Then $\mathcal{Z}\text{-}\lim_n x_n=\eta$ if and only if $\lim_n x_n=\eta$. 
\end{thm}

Fridy also proved in the same work that the growth rate $O(1/n)$ is \textquotedblleft best possible\textquotedblright\, in the following sense, see \cite[Theorem 4]{MR816582}:
\begin{thm}\label{thm:conversefridy}
Let $r$ be a decreasing positive real sequence such that $r_n \neq O(1/n)$. 
Then there exists a nonconvergent real sequence $x$ such that $|\Delta x_n|\le r_n$ for all $n \in \omega$ and $\mathcal{Z}\text{-}\lim_n x_n=0$. 
\end{thm}
Analogous Tauberian conditions related to statistical convergence have been obtained by Maddox in \cite{MR938459, MR1002541}.  
Here we study Tauberian conditions for arbitrary ideals and
for the preservation of ideal cluster and limit points.

In Section \ref{sec:mainresults} we state our main results. Rather informally, we will show the following (here, we can think $(r_n)_{n \in \omega}$ as a positive null sequence with divergent sum): 
\begin{enumerate}[label={\rm (\roman{*})}]
\item (Theorem \ref{equivalencefridyclaim}) $|\Delta x_n|=O(r_n)$ is a Tauberian condition for the equivalence between ordinary convergence and $\mathcal{I}$-convergence if and only if every accumulation point is an $\mathcal{I}$-cluster point (for an arbitrary ideal $\mathcal{I}$); an analogue is given for $\mathcal{I}$-limit points; 
\item (Theorem \ref{thm:charactbaire}) A Tauberian condition $|\Delta x_n|=O(r_n)$ for the equivalence between ordinary convergence and $\mathcal{I}$-convergence exists if and only if the ideal $\mathcal{I}$ has the Baire property; 
\item (Theorem \ref{thm:mainanalyticPideal}) Within the family of analytic $P$-ideals, we characterize when $|\Delta x_n|=O(1/n)$ is a Tauberian condition for the equivalence between ordinary convergence and ideal convergence; 
\item A notion of \textquotedblleft optimal Tauberian rate\textquotedblright\, is studied in Section \ref{sec:optimal}. In Proposition \ref{prop:ratesummableo}, we show that, for each suitable $(r_n)$, there exists a summable ideal $\mathcal{I}$ such that $|\Delta x_n|=o(r_n)$ is an optimal Tauberian rate for both $\Lambda_x(\mathcal{I})=\mathrm{L}_x$ and $\Gamma_x(\mathcal{I})=\mathrm{L}_x$. Lastly, we show that Fridy's Theorem \ref{thm:conversefridy} can be extended to all ideals $\mathcal{I}$ with $\mathcal{I}_{1/n}\subseteq \mathcal{I}\subseteq \mathcal{Z}$ (Theorem \ref{thm:optimalsummableideal}). 
\end{enumerate}
Some preliminaries are given in Section \ref{sec:preliminaries}, while the proofs of our main results are given in Section \ref{sec:proofs}.







\section{Main results} \label{sec:mainresults}

\subsection{Equivalences of claims.}\label{subsec:equivalences}
We start by showing that a Tauberian condition on the growth rate of $\Delta x$ which ensures the equivalence between $\mathcal{I}$-convergence and ordinary convergence is equivalent to the requirement $\Gamma_x(\mathcal{I})=\mathrm{L}_x$.

Recall that an ideal $\mathcal{I}$ is a $P^{+}$\emph{-ideal} if for every decreasing sequence $(S_n)_{n \in \omega}$ of sets in $\mathcal{I}^+$, there exists $S \in \mathcal{I}^+$ such that $S\setminus S_n$ is finite for all integers $n \in \omega$. Examples of $P^{+}$-ideals include all $F_\sigma$ ideals (and, more generally, all $G_{\delta\sigma}$ ideals) by \cite[Proposition 3.2(iii)]{FKL2025}. 
It is consistent (e.g.,~assuming the Continuum Hypothesis) that there exists a maximal ideal which is $P^+$, see \cite{MR0080902}. However, it is also consistent that there are no maximal ideals which are $P^+$, see \cite{MR0728877}. 
On the other hand, it is remarkable that the ideal $\mathcal{Z}$ defined in \eqref{eq:definitionZ} is not a $P^{+}$-ideal, see \cite[Example 3.3(ii)]{FKL2025}.

\begin{thm}\label{equivalencefridyclaim}
Let $r$ be a positive real sequence such that $\lim_n r_n=0$ and $\sum_n r_n=\infty$. Let also $\mathcal{I}$ be an ideal on $\omega$. Consider the following statements\textup{:} 
\begin{enumerate}[label={\rm (\textsc{A}\arabic{*})}]
\item \label{item:A1} if $x$ is a real sequence such that 
$|\Delta x_n|=O(r_n)$, 
then $\Lambda_x(\mathcal{I})=\mathrm{L}_x$\textup{;}
\item \label{item:A2} if $x$ is a real sequence such that 
$|\Delta x_n|=O(r_n)$, 
then $\Gamma_x(\mathcal{I})=\mathrm{L}_x$\textup{;}
\item \label{item:A3} if $x$ is a real sequence such that 
$|\Delta x_n|=O(r_n)$ and $\mathcal{I}\text{-}\lim x=\eta$, 
then $\lim x=\eta$\textup{.}
\end{enumerate}
Then \ref{item:A1} $\implies$ \ref{item:A2} $\Longleftrightarrow$ \ref{item:A3}. 

If, in addition, $\mathcal{I}$ is a $P^+$-ideal, then \ref{item:A1} $\Longleftrightarrow$ \ref{item:A2}. 
\end{thm}


Observe that the case where $\lim_nr_n=0$ and $\sum_n r_n<\infty$ is not interesting. In fact, in such case, if $x$ is a real sequence such that $|\Delta x_n|=O(r_n)$ then $x$ would be Cauchy since for all $n,m \in \omega$ with $n\le m$ it holds $|x_n-x_m| \le \sum_{k=n}^{m-1}|\Delta x_k| \le \sum_{k\ge n}|\Delta x_k|$. Hence $x$ would be convergent to some $\eta \in \mathbb{R}$, so that $\Lambda_x(\mathcal{I})=\Gamma_x(\mathcal{I})=\mathrm{L}_x=\{\eta\}$. 

\begin{rmk}\label{rmk:Imaximal}
%
It is worth noting that, if $\mathcal{I}$ is a maximal ideal on $\omega$ then the equivalence \ref{item:A1} $\Longleftrightarrow$ \ref{item:A2} holds:  
to this aim, let $(I_k)_{k \in \omega}$ be a partition of $\omega$ into consecutive finite intervals such that $\sum_{n \in I_k}r_n>2+3\max\{r_n: n\in \omega\}$ for all $k \in \omega$ (which is possible since $\lim_n r_n=0$ and $\sum_n r_n=\infty$). 
It follows by Lemma \ref{lem:technicalbump} below that there exists a real sequence $x$ such that, for each $k \in \omega$, the following hold: $x_{\min I_k}=x_{\max I_k}=0$, $x_n=1$ for some $n \in I_k$, and $|\Delta x_n| \le r_n$ for all $n \in I_k$. It follows by construction and the maximality of $\mathcal{I}$ that 
$$
\mathrm{L}_x=[0,1], \quad \Lambda_x(\mathcal{I})\subseteq \Gamma_x(\mathcal{I}), 
 \quad \text{ and }\quad 
|\Gamma_x(\mathcal{I})|=1.
$$
Therefore both items \ref{item:A1} and \ref{item:A2} fail. 
\end{rmk}

We observe also that the equality $\Gamma_x(\mathcal{I})=\mathrm{L}_x$ in item \ref{item:A2} \emph{implies} the equivalence between $\mathcal{I}$-convergence and ordinary convergence in item \ref{item:A3} in the case where $x$ is bounded, see e.g. \cite[Corollary 3.4]{MR3920799}. However, if $x$ is unbounded, the latter implication relies on the additional growth condition of the sequence $\Delta x$. 

\begin{rmk}
Given an ideal $\mathcal{I}$ on $\omega$, let $c(\mathcal{I})$ be the vector space of $\mathcal{I}$-convergent real sequences. Observe that the operator $c(\mathcal{I})\to \mathbb{R}$ defined by $x\mapsto \mathcal{I}\text{-}\lim x$ is linear. In addition, for each real sequence $x$ and nonzero constant $c \in \mathbb{R}$, we have $\Gamma_{(x_n/c)}(\mathcal{I})=\{\eta/c: \eta \in \Gamma_x(\mathcal{I})\}$ and, similarly, $\Lambda_{(x_n/c)}(\mathcal{I})=\{\eta/c: \eta \in \Lambda_x(\mathcal{I})\}$, cf. \cite[Proposition 3.2]{MR4428911}. It easily follows that items \ref{item:A1}, \ref{item:A2}, and \ref{item:A3} are also equivalent to its corresponding primed statement: 
\begin{enumerate}[label={\rm (\textsc{A}\arabic{*}$^\prime$)}]
\item \label{item:A1prime} if $x$ is a real sequence such that 
$|\Delta x_n|\le r_n$ for all $n$, 
then $\Lambda_x(\mathcal{I})=\mathrm{L}_x$\textup{;}
\item \label{item:A2prime} if $x$ is a real sequence such that 
$|\Delta x_n|\le r_n$ for all $n$, 
then $\Gamma_x(\mathcal{I})=\mathrm{L}_x$\textup{;}
\item \label{item:A3prime} if $x$ is a real sequence such that 
$|\Delta x_n|\le r_n$ for all $n$,  and $\mathcal{I}\text{-}\lim x=\eta$, 
then $\lim x=\eta$\textup{.}
\end{enumerate}
\end{rmk}

In turn, Remark \ref{rmk:Imaximal} implies that the equivalence \ref{item:A1} $\Longleftrightarrow$ \ref{item:A2} does not hold only for $P^+$-ideals: for instance, let $\mathcal{I}$ be a maximal ideal which is not $P^+$ (which exists in ZFC: for instance, pick any maximal ideal extending $\mathrm{Fin}^2$, see e.g. \cite[Example 3.4(ii)]{FKL2025}). 
Additional examples of non $P^+$-ideal which satisfy the equivalence  \ref{item:A1} $\Longleftrightarrow$ \ref{item:A2} can be obtained by Theorem \ref{thm:connectedBanach} below. 

We will show in Example \ref{ex:nonequivalence} below that the equivalence \ref{item:A1} $\Longleftrightarrow$ \ref{item:A2} in Theorem \ref{equivalencefridyclaim} fails for the ideal $\mathcal{Z}$. 
We leave as open question for the interested reader to characterize, at a given rate $r$, the class of ideals $\mathcal{I}$ for which the equivalence \ref{item:A1} $\Longleftrightarrow$ \ref{item:A2} holds. 
%

\subsection{Baire ideals and arbitrary growth} 
In the next result, we show that the Tauberian condition \ref{item:A2} for ideal cluster points (or, equivalently, \ref{item:A3} for ideal convergence) holds for some sequence $r$ if and only if the ideal $\mathcal{I}$ has the Baire property. In the proof, we will use Talagrand's characterization, namely, an ideal $\mathcal{I}$ on $\omega$ has the Baire property if and only if there exists a partition $(A_j: j \in \omega)$ of $\omega$ such that each $A_j$ is a finite interval of nonnegative integers and $S\notin \mathcal{I}$ whenever $A_j \subseteq S$ for infinitely many $j \in \omega$; see \cite[Theorem 2.1]{MR579439} and also \cite[Theorem 4.1.2]{MR1350295}. Additional characterizations of the Baire property for ideals can be found in \cite{MR4566746, MR4165727, MR4915070}. 
\begin{thm}
\label{thm:charactbaire}
Let $\mathcal{I}$ be an ideal on $\omega$. Then the following are equivalent\textup{:}
\begin{enumerate}[label={\rm (\textsc{B}\arabic{*})}]
\item \label{item:BP1} $\mathcal{I}$ has the Baire property\textup{;}
\item \label{item:BP2A} there exists a positive real sequence $r$ for which $\lim_n r_n=0$, $\sum_n r_n=\infty$, and for which, if $x$ is a real sequence such that 
$|\Delta x_n|=O(r_n)$, 
then $\Lambda_x(\mathcal{I})=\mathrm{L}_x$\textup{;}
\item \label{item:BP2} there exists a positive real sequence $r$ for which $\lim_n r_n=0$, $\sum_n r_n=\infty$, and for which, if $x$ is a real sequence such that 
$|\Delta x_n|=O(r_n)$, 
then $\Gamma_x(\mathcal{I})=\mathrm{L}_x$\textup{;}
\item \label{item:BP3} there exists a positive real sequence $r$ for which $\lim_n r_n=0$, $\sum_n r_n=\infty$, and for which, if $x$ is a real sequence such that 
$|\Delta x_n|=O(r_n)$ and $\mathcal{I}\text{-}\lim x=\eta$, 
then $\lim x=\eta$\textup{.}
\end{enumerate}
\end{thm}

\begin{rmk}
\label{rmk:choicern}
As follows from \eqref{eq:kjdshgsd} in the proof of the implication \ref{item:BP1} $\implies$ \ref{item:BP2A}, if $\mathcal{I}$ is an ideal with the Baire property and $(A_j: j \in \omega)$ is a sequence witnessing Talagrand's characterization, then it is enough to choose a suitable $r$ with 
$$
r_n = o\left(\frac{1}{|A_{\iota(n)}|}\right), 
$$
where $\iota: \omega\to \omega$ is the (weakly) increasing function such that $n \in A_{\iota(n)}$ for all $n \in \omega$. In addition, the same sequence $r$ works also for items \ref{item:BP2} and \ref{item:BP3} by Theorem \ref{equivalencefridyclaim}. 
\end{rmk}

\begin{thm}
\label{thm:conjectureweak}
Let $x$ be a real sequence such that $|\Delta x_n|=o(1/n)$, and fix $\eta \in \mathbb{R}$. 
Then $\Lambda_x(\mathcal{Z})
=\mathrm{L}_x$, and $\mathcal{Z}\text{-}\lim_n x_n=\eta$ if and only if $\lim_n x_n=\eta$. 
\end{thm}
It is remarkable that the Tauberian condition $|\Delta x_n|=o(1/n)$ is \textquotedblleft optimal\textquotedblright\,(cf. Section \ref{sec:optimal} below) for the equality $\Lambda_x(\mathcal{Z})=\mathrm{L}_x$. In fact, as it follows from Example \ref{ex:nonequivalence} below, there exists a real sequence $x$ such that $|\Delta x_n|=O(1/n)$ and $\Lambda_x(\mathcal{Z})\neq \mathrm{L}_x$. 

In addition, since $\Lambda_x(\mathcal{Z})\subseteq \Gamma_x(\mathcal{Z})\subseteq \mathrm{L}_x$, it follows that if $|\Delta x_n|=o(1/n)$, then also $\Gamma_x(\mathcal{Z})= \mathrm{L}_x$. Not only this is weaker than Fridy's Theorem \ref{thm:fridythm3}, but $o(1/n)$ is \emph{not} an optimal rate for the latter equality: in fact, as  it follows from Theorem \ref{equivalencefridyclaim} and Fridy's Theorem \ref{thm:fridythm3} and Theorem \ref{thm:conversefridy}, an optimal growth for $\Gamma_x(\mathcal{Z})= \mathrm{L}_x$ is $O(1/n)$ among nonincreasing admissible rates. With a different method, we will recover the rate $O(1/n)$ in Corollary \ref{cor:originalconjecture} below. 

Lastly, even if the second part of Theorem \ref{thm:conjectureweak} is weaker than Fridy's Theorem \ref{thm:fridythm3} (since it requires the stronger condition $o(1/n)$), one has to take in mind that, right now, we have been using only the fact that the ideal $\mathcal{Z}$ has the Baire property.

%

\begin{rmk}\label{rmk:o1overnBaire}
With the same argument as in Theorem \ref{thm:conjectureweak}, we have that the analogous conclusion holds replacing $\mathcal{Z}$ with an ideal $\mathcal{I}$ for which there exists a real $\alpha>1$ such that $\bigcup_{j \in J} (\omega \cap [\alpha^j,\alpha^{j+1}))\in \mathcal{I}^+$ for all infinite $J\subseteq \omega$.  For instance, this applies if $\mathcal{I}$ is the summable ideal $\{S\subseteq \omega: \sum_{n \in S\setminus \{0\}}1/n <\infty\}$. 
\end{rmk}






Given an arbitrary positive real sequence $r$  as in the statement of Theorem \ref{equivalencefridyclaim}, one might also ask whether there exists an ideal $\mathcal{I}\neq \mathrm{Fin}$, possibly with sufficiently low topological complexity, which satisfies the (strongest) item \ref{item:A1}. We answer this question in the affirmative within the class of ideals with lowest topological complexity (recall that, as a consequence of Baire's category theorem, an ideal on $\omega$ cannot be $G_\delta$, see e.g. \cite{alcantara-phd-thesis}).
\begin{prop}\label{prop:Baire}
Let $r$ be a positive real sequence such that $\lim_n r_n=0$. 
Then there exists an $F_\sigma$ ideal $\mathcal{I}$ on $\omega$ which is not a $P$-ideal and such that, if $x$ is a real sequence such that $|\Delta x_n|=O(r_n)$, then $\Lambda_x(\mathcal{I})=
\mathrm{L}_x$\textup{.}
\end{prop}
Since $\mathcal{I}$ is not a $P$-ideal then necessarily $\mathcal{I}$ is not summable (in particular, $\mathcal{I}\neq \mathrm{Fin}$). Hence, Proposition \ref{prop:Baire} tells us that there are nontrivial cases with arbitrarily large Tauberian conditions. 


\subsection{Analytic P-ideals and growth $O(1/n)$}\label{sec:analyticP} 
Let us recall that an ideal $\mathcal{I}$ on $\omega$ is a $P$\emph{-ideal} 
if it is $\sigma$-directed modulo finite sets, that is, for every sequence $(S_n)$ of sets in $\mathcal{I}$ there exists $S \in \mathcal{I}$ such that $S_n\setminus S$ is finite for all integers $n\in \omega$. Examples of $P$-ideals which are also analytic (as subsets of the Cantor space) include: $\mathrm{Fin}$, $\mathcal{Z}$, all summable ideals (that is, ideals of the type 
\begin{equation}\label{eq:summable}
\mathcal{I}:=\left\{S\subseteq \omega: \sum_{n \in S}a_n<\infty\right\}
\end{equation}
for some nonnegative real sequence $(a_n)_{n \in \omega}$ with $\sum_n a_n=\infty$), ideals which are generated by nonnegative regular matrices \cite{Filipow18}, 
the Fubini product $\emptyset \times \mathrm{Fin}$, which can be defined as 
$
\{S\subseteq \omega: \forall n \in \omega, S \cap I_n \in \mathrm{Fin}\},
$ 
where $(I_n)$ is a given partition of $\omega$ into infinite sets, 
Erd{\H o}s--Ulam ideals introduced by Just and Krawczyk in \cite{MR748847}, 
certain ideals used by Louveau and Veli\u{c}kovi\'{c} \cite{Louveau1994}, 
and, more generally, density-like ideals and generalized density ideals \cite{MR3436368, MR4404626}. 
Additional pathological examples can be found in \cite{MR0593624}. 
It has been suggested in \cite{MR4124855, MR3436368} that the theory of analytic $P$-ideals may have some relevant yet unexploited potential for the study of the geometry of Banach spaces.

Hereafter, a function $\varphi: \mathcal{P}(\omega) \to [0,\infty]$ is a  \emph{lower semicontinuous submeasure} (in short, \emph{lscsm}) if the following properties hold: 
\begin{enumerate}[label={\rm (\roman{*})}]
\item $\varphi(\emptyset)=0$, 
\item $\varphi(A)\le \varphi(B)$ for all $A\subseteq B\subseteq \omega$,
\item $\varphi(A\cup B)\le \varphi(A)+\varphi(B)$ for all $A,B\subseteq \omega$, 
\item $\varphi(F)<\infty$ for all $F \in \mathrm{Fin}$, and 
\item $\varphi(A)=\sup\{\varphi(A\cap F): F \in \mathrm{Fin}\}$ for all $A\subseteq \omega$. 
\end{enumerate}
By a classical result of Solecki, see \cite[Theorem 3.1]{MR1708146}, an ideal $\mathcal{I}$ on $\omega$ is an analytic $P$-ideal if and only if there exists a lscsm 
$\varphi: \mathcal{P}(\omega) \to [0,\infty]$ 
such that
\begin{equation}\label{eq:characterizationanalPideal}
\mathcal{I}=\mathrm{Exh}(\varphi):=\{S\subseteq \omega: \|S\|_\varphi=0\} 
\,\,\,\text{ and }\,\,\,
0<\|\omega\|_\varphi \le \varphi(\omega)<\infty,
\end{equation}
where 
$$
\|S\|_\varphi:=\inf\{\,\varphi(S\setminus F): F \in \mathrm{Fin}\}
$$
for all $S\subseteq \omega$. In particular, every analytic $P$-ideal is necessarily $F_{\sigma\delta}$. 
Note that $\|\cdot\|_\varphi$ is a monotone subadditive function which is invariant modulo finite sets. Also, the choice of $\varphi$ in \eqref{eq:characterizationanalPideal} is not necessarily unique (even up to constants): for instance $\mathcal{Z}=\mathrm{Exh}(\varphi)=\mathrm{Exh}(\psi)$, where $\varphi$ and $\psi$ are lscsms defined by $\varphi(S)=\sup_n |S\cap [0,n]|/(n+1)$ and $\psi(S)=\sup_n |S\cap [2^n,2^{n+1})|/2^n$ for each $S\subseteq \omega$. 

We say also that a lscsm $\varphi: \mathcal{P}(\omega)\to [0,\infty]$ is \emph{smooth} if: 
\begin{equation}\label{eq:definitionsmooth}
\forall \kappa>1, \quad 
\liminf_{n\to \infty}\varphi([n, \kappa n])>0.
\end{equation}
As it follows from Theorem \ref{thm:mainanalyticPideal}, if $\varphi$ and $\varphi^\prime$ are lscsms such that $\mathrm{Exh}(\varphi)=\mathrm{Exh}(\varphi^\prime)$ then $\varphi$ is smooth if and only if $\varphi^\prime$ is smooth as well. In addition, if $\mathcal{I}$ is represented as in \eqref{eq:characterizationanalPideal}, we show that the notion of smoothness is equivalent to the fact that every accumulation point of $x$ is an $\mathcal{I}$-cluster point, for every sequence $x$ with growth rate $|\Delta x_n|=O(1/n)$. 
\begin{thm}\label{thm:mainanalyticPideal}
Let $\mathcal{I}$ be an analytic $P$-ideal on $\omega$. 
Then the following are equivalent\textup{:}
\begin{enumerate}[label={\rm (\textsc{P}\arabic{*})}]
\item \label{item:Analytic1} if $x$ is a real sequence such that 
$|\Delta x_n|=O(1/n)$, 
then $\Gamma_x(\mathcal{I})=\mathrm{L}_x$\textup{;}
\item \label{item:Analytic1bis} if $x$ is a real sequence such that 
$|\Delta x_n|=O(1/n)$ 
 and $\mathcal{I}\text{-}\lim_n x_n=\eta$, then $\lim_n x_n=\eta$\textup{;} 
\item \label{item:Analytic2} if $\varphi$ is a lscsm such that $\mathcal{I}=\mathrm{Exh}(\varphi)$ then it is smooth\textup{;}
\item \label{item:Analytic3} there exists a smooth lscsm $\varphi$ such that $\mathcal{I}=\mathrm{Exh}(\varphi)$\textup{.}
\end{enumerate}
\end{thm}

As an application of Theorem \ref{thm:mainanalyticPideal}, we recover Fridy's Theorem \ref{thm:fridythm3} (by showing that there exists at least one smooth lscsm $\varphi$ which generates $\mathcal{Z}$): 
\begin{cor}\label{cor:originalconjecture}
Let $\mathcal{I}$ be an ideal on $\omega$ such that $\mathcal{I}\subseteq \mathcal{Z}$. 
Also, let $x$ be a real sequence such that $|\Delta x_n|=O(1/n)$ 
and fix $\eta \in \mathbb{R}$. 

Then $\Gamma_x(\mathcal{I})=\mathrm{L}_x$, and $\mathcal{I}\text{-}\lim_n x_n=\eta$ if and only if $\lim_n x_n=\eta$. 
\end{cor}

As another application, given an ideal $\mathcal{I}$ on $\omega$, let us recall that a real sequence $x=(x_n)$ is $\mathcal{I}$\emph{-statistically convergent} to $\eta$ if 
$$
\forall \varepsilon>0, \quad 
\mathcal{I}\text{-}\lim_{n\to \infty} \lambda_n\left(\left\{k \in \omega: |x_k-\eta|>\varepsilon \right\}\right)=0,
$$
where $\lambda_n(S):=|S\cap [0,n]|/(n+1)$ for all $n \in \omega$ and $S\subseteq \omega$, see \cite[Definition 2.2]{MR4052262}. We obtain that $\mathcal{Z}$-statistical convergence is nothing but statistical convergence (that is, $\mathcal{Z}$-convergence), recovering \cite[Corollary 2.8]{MR4052262}:
\begin{cor}
\label{cor:recovercorollary}
Let $\mathcal{I}$ be an ideal on $\omega$ such that $\mathcal{I}\subseteq \mathcal{Z}$. Then $\mathcal{I}$-statistical convergence coincides with statistical convergence \textup{(}i.e., $\mathcal{Z}$-convergence\textup{)}. 
\end{cor}

Unlike the case of the ideal $\mathcal{Z}$ considered in Corollary \ref{cor:originalconjecture}, it is remarkable that not all analytic $P$-ideals are generated by some smooth lscsm. In fact, we show below that there exists a summable ideal as in \eqref{eq:summable} (which is an $F_\sigma$ $P$-ideal, hence also analytic) that satisfies $\Gamma_x(\mathcal{I})=\mathrm{L}_x$ for all sequences with growth rate $o(1/n)$ but not necessarily for the ones with growth rate $O(1/n)$. In this sense, the rate $o(1/n)$ will be the \textquotedblleft optimal\textquotedblright\, Tauberian rate for such an ideal (for more on optimal rates, see \S \ref{sec:optimal} below).
\begin{prop}\label{prop:summableborderline}
There exists a summable ideal $\mathcal{I}$ on $\omega$ with the following properties\textup{:}
\begin{enumerate}[label={\rm (\roman{*})}]
\item \label{item:example1fixed} if $x$ is a real sequence such that 
$|\Delta x_n|=o(1/n)$, then $
\Gamma_x(\mathcal{I})= \mathrm{L}_x$\textup{;} 
\item \label{item:example2fixed} there is no smooth lscsm $\varphi$ such that $\mathrm{Exh}(\varphi)=\mathcal{I}$\textup{;}
\item \label{item:example3fixed} there exists a real sequence $x$ such that 
$|\Delta x_n|=O(1/n)$ and $\Gamma_x(\mathcal{I})\neq \mathrm{L}_x$\textup{.} 
\end{enumerate}
\end{prop}

We conclude this section by 
proving that the equivalence \ref{item:A1} $\Longleftrightarrow$ \ref{item:A2} in Theorem \ref{equivalencefridyclaim} fails for the ideal $\mathcal{Z}$, namely, there exists a sequence $x$ with $|\Delta x_n|=O(1/n)$ such that 
$$
\Gamma_x(\mathcal{Z})=\mathrm{L}_x=[0,1]
\quad \text{ and }\quad 
\Lambda_x(\mathcal{Z})=\emptyset.
$$
\begin{example}\label{ex:nonequivalence}
Define $A_0:=\{0,1\}$ and for each nonzero integer $k\in \omega$ set $A_k:=\omega \cap [2^k,2^{k+1})$. 
Now, let $x=(x_n)_{n \in \omega}$ be the real sequence defined by $x_n:=0$ if $n \in A_0$ and, for each nonzero integer $k \in \omega$, set 
$$
\forall i \in \omega \cap [0,2^{k-1}), \quad 
x_{2^k+i}:=x_{2^{k+1}-1-i}:=1-\frac{i}{2^{k-1}}.
$$
For instance, the restriction of $x$ to the interval $A_3$, that is, $(x_8,x_9,\ldots,x_{15})$ is defined as $(1, \nicefrac{3}{4},\nicefrac{1}{2}, \nicefrac{1}{4}, \nicefrac{1}{4}, \nicefrac{1}{2}, \nicefrac{3}{4}, 1)$. It is easy to see that 
$$
\mathrm{L}_x=[0,1]
\quad \text{ and }\quad 
|\Delta x_n|\le \frac{1}{2^{\lfloor \log_2(n)\rfloor-1}}=O\left(\frac{1}{n}\right) \text{ for }n\ge 1.
$$
It follows by Corollary \ref{cor:originalconjecture} that $\Gamma_x(\mathcal{Z})=\mathrm{L}_x$. 

On the other hand, fix $\eta\in\mathrm{L}_x=[0,1]$ and pick an
infinite $K\subseteq\omega$ such that $(x_n:n\in K)$ is convergent
to $\eta$. Then for each $j\in\omega$ there exists a finite set
$F_j\in\mathrm{Fin}$ such that $K\subseteq T_j\cup F_j$, where
$T_j:=\{i\in\omega:|x_i-\eta|\le 2^{-j}\}$.
Since the values on each half of $A_k$ are spaced $2^{1-k}$
apart, we have $|T_j\cap A_k|\le 2(2^{k-j}+1)$ for all $k\ge1$.
Recalling that $\mathcal{Z}=\mathrm{Exh}(\psi)$, where
$\psi:\mathcal{P}(\omega)\to[0,\infty]$ is the lscsm defined by
$\psi(S)=\sup_n |S\cap[2^n,2^{n+1})|/2^n$ for each
$S\subseteq\omega$, it follows that
\begin{displaymath}
\begin{split}
\|K\|_\psi
\le \|T_j\|_\psi+\|F_j\|_\psi
=\limsup_{k\to\infty}\frac{|T_j\cap A_k|}{|A_k|}\le \limsup_{k\to\infty}\frac{2(2^{k-j}+1)}{2^k}
=\frac{1}{2^{j-1}}
\end{split}
\end{displaymath}
for all nonzero $j\in\omega$.
Thus $K\in\mathcal{Z}$, so that
$\eta\notin\Lambda_x(\mathcal{Z})$.
Since $\eta\in\mathrm{L}_x$ was arbitrary and
$\Lambda_x(\mathcal{Z})\subseteq\mathrm{L}_x$, we conclude that
$\Lambda_x(\mathcal{Z})=\emptyset$. In particular, $\Lambda_x(\mathcal{Z})\neq \mathrm{L}_x$. 
\end{example}

\subsection{Optimal Tauberian rates} \label{sec:optimal} 
Let us start with the definition which will allow us to speak of   \textquotedblleft optimal Tauberian rates.\textquotedblright\, To  this aim and to ease the notation, let us write 
$$
\mathcal{R}:=\left\{(r_n)_{n \in \omega} \in (0,\infty)^\omega: \lim_{n\to \infty} r_n=0 \,\, \text{ and }\,\,\sum_{n\in\omega} r_n=\infty\right\}.
$$
Let also $\mathcal{R}^{\downarrow}$ be the subset of $\mathcal{R}$ containing only its (weakly) decreasing sequences $(r_n)_{n \in \omega}$. 

\begin{defi}\label{defi:optimality}
Let $\mathscr{P}$ be a statement about a real sequence $x$, which can be either true or false. 
Fix a sequence $(r_n)_{n \in \omega} \in (0,\infty)^\omega$. 
 %
We say that $o(r_n)$ \emph{is an optimal Tauberian rate for} $\mathscr{P}$ 
if:  
\begin{enumerate}[label={\rm 
(\roman{*})}]
\item $\mathscr{P}$ holds whenever $x$ is a real sequence such that $|\Delta x_n|=o(r_n)$; 
\item there exists a sequence $x$ such that $|\Delta x_n|= O(r_n)$ and $\mathscr{P}$ fails. 
\end{enumerate}
Also, let $\mathcal{F}\subseteq (0,\infty)^\omega$ be a nonempty family of sequences. We say that $O(r_n)$ \emph{is an optimal Tauberian rate for} $\mathscr{P}$ \emph{within} $\mathcal{F}$ 
if:  
\begin{enumerate}[label={\rm 
(\roman{*})}]
\item $\mathscr{P}$ holds whenever $x$ is a real sequence such that $|\Delta x_n|=O(r_n)$; 
\item for all $r^\prime \in \mathcal{F}$ with $r^\prime_n \neq O(r_n)$, there exists a sequence $x$ with $|\Delta x_n|= O(r^\prime_n)$ such that $\mathscr{P}$ fails. 
\end{enumerate}
\end{defi}

With the above definition at hand, by Proposition \ref{prop:summableborderline} there exists a summable ideal $\mathcal{I}$ such that 
$
o\left(1/n\right) \text{ is an optimal rate for }\Gamma_x(\mathcal{I})=\mathrm{L}_x.
$ 
At the same time, recall that summable ideals are $F_\sigma$, and that all $F_\sigma$ ideals (and, more generally, all $G_{\delta\sigma}$) are $P^+$ by \cite[Proposition 3.2(iii)]{FKL2025}. By the equivalence \ref{item:A1} $\Longleftrightarrow$ \ref{item:A2} in Theorem \ref{equivalencefridyclaim} we obtain, for this summable ideal $\mathcal{I}$, that: 
$$
o\left(1/n\right) \text{ is an optimal rate for both }\Lambda_x(\mathcal{I})=\mathrm{L}_x 
\text{ and }\Gamma_x(\mathcal{I})=\mathrm{L}_x.
$$

In addition, taking into account the equivalence \ref{item:A2} $\Longleftrightarrow$ \ref{item:A3} in Theorem \ref{equivalencefridyclaim}, Fridy's Theorem \ref{thm:fridythm3} and Theorem \ref{thm:conversefridy} can be simply rewritten as: 
$$
O\left(1/n\right) \text{ is an optimal rate for }\Gamma_x(\mathcal{Z})=\mathrm{L}_x \text{ within }\mathcal{R}^{\downarrow}.
$$
As we already observed, optimal rates for ideal limit points can be different. In fact, as it follows by Theorem \ref{thm:conjectureweak} and Example \ref{ex:nonequivalence}, we have: 
$$
o\left(1/n\right) \text{ is an optimal rate for }\Lambda_x(\mathcal{Z})=\mathrm{L}_x. 
$$


\medskip

It is remarkable that Fridy's Theorem \ref{thm:conversefridy} relies heavily on the fact that the sequence $r$ is decreasing. In fact, in the following example we show that: 
$$
O(1/n) \text{ is } not \text{ an optimal rate for } \Gamma_x(\mathcal{Z})=\mathrm{L}_x 
\text{ within }\mathcal{R}.
$$
\begin{example}\label{example:aboutfridyconverse}
Define 
the positive real sequence $r^\prime \in \mathcal{R}$ by 
\begin{displaymath}
r^\prime_n:=
\begin{cases}
\,
1/(k+1) 
\,\,\,&\text{ if }n=2^k \text{ for some }k \in \omega\\
\,
1/3^n
 & \text{ otherwise}. 
\end{cases}
\end{displaymath}
It is immediate that $r^\prime$ is not decreasing, so that $r^\prime\notin \mathcal{R}^\downarrow$, and that $r^\prime_n \neq O(1/n)$. Define the intervals $A_0:=\{0,1\}$ and $A_k:=[2^k,2^{k+1})$ for all $k \ge 1$, and observe that $(A_k)_{k \in \omega}$ is a witnessing sequence for the Baire property of 
$\mathcal{Z}$. In addition, 
$$
\sum_{n \in A_k}r^\prime_n \le \frac{1}{k+1}+\frac{|A_k|}{3^k}=o(1).
$$
It follows by Proposition \ref{prop:basicestimateLambdaL} that, if $x$ is a real sequence with $|\Delta x_n|=O(r^\prime_n)$, then $\Lambda_x(\mathcal{Z})=\mathrm{L}_x$. In particular, $O(r^\prime_n)$ is a Tauberian condition for $\Gamma_x(\mathcal{Z})=\mathrm{L}_x$. 
\end{example}

A somehow related distinction between $o$- and $O$-Tauberian conditions already appears in the theory of matrix summability, see Lorentz and Zeller \cite{MR346370}.

Now, we show that every rate $o(r_n)$ can be realized as an optimal rate of some summable ideal as in \eqref{eq:summable}.
\begin{prop}\label{prop:ratesummableo}
Fix $r \in \mathcal{R}$. Then there exists a summable ideal $\mathcal{I}$ on $\omega$ such that $o(r_n)$ is an optimal rate for both $\Lambda_x(\mathcal{I})=\mathrm{L}_x$ and $\Gamma_x(\mathcal{I})=\mathrm{L}_x$. 
\end{prop}

Summable ideals also admit optimal rates of the type $O(1/n)$ within $\mathcal{R}^{\downarrow}$, cf. Remark \ref{rmk:oandO} below. To this aim, consider the summable ideal
$$
\mathcal{I}_{1/n}:=\left\{S\subseteq \omega: \sum_{n \in S\setminus \{0\}}1/n<\infty\right\}. 
$$  
It is folklore that $\mathcal{I}_{1/n}\subseteq \mathcal{Z}$ (and the inclusion is strict: for instance, the set of primes belongs to $\mathcal{Z}\setminus \mathcal{I}_{1/n}$). Then, taking into account Theorem \ref{equivalencefridyclaim}, we have the following result, which generalizes both Fridy's Theorem \ref{thm:fridythm3} and Theorem \ref{thm:conversefridy}. 
\begin{thm}\label{thm:optimalsummableideal}
Let $\mathcal{I}$ be an ideal on $\omega$ such that $\mathcal{I}_{1/n}\subseteq \mathcal{I}\subseteq \mathcal{Z}$. Then 
$O(1/n)$ is an optimal rate for $\Gamma_x(\mathcal{I})=\mathrm{L}_x$ within $\mathcal{R}^{\downarrow}$. 
\end{thm}

\begin{rmk}\label{rmk:oandO}
The two notions in Definition \ref{defi:optimality} are \emph{not}
mutually exclusive. Indeed, let $J_k:=[4^k,4^k+k)\cap\omega$
for all $k\ge8$, and define $r\in\mathcal{R}$ by
$$
r_n:=
\begin{cases}
\,1/k \,\,\,\,\,& \text{if }n\in J_k\text{ for some }k\ge8,\\
\,1/(n+1) & \text{otherwise}.
\end{cases}
$$
The same partition $(A_k)$ of Example
\ref{example:aboutfridyconverse} witnesses the Baire property
of $\mathcal{I}_{1/n}$ and satisfies
$\sup_k\sum_{n\in A_k}r_n<\infty$. Since $\mathcal{I}_{1/n}$ is an $F_\sigma$-ideal, then $\Lambda_x(\mathcal{I}_{1/n})=\Gamma_x(\mathcal{I}_{1/n})$ for every sequence $x$. 
Hence Proposition \ref{prop:basic2estimateLambdaLwitho} below yields
$\Gamma_x(\mathcal{I}_{1/n})=\mathrm{L}_x$ whenever
$|\Delta x_n|=o(r_n)$.
On the other hand,
$$
\sum_{k\ge8}\sum_{n\in J_k}\frac1n
\le\sum_{k\ge8}\frac{k}{4^k}<\infty.
$$
Applying Lemma \ref{lem:technicalbump} below on each $J_k$ with $\eta:=0$ and $\delta:=1/4$, and setting $x_n:=0$ elsewhere,
we obtain $|\Delta x_n|\le r_n$,
$\mathcal{I}_{1/n}\text{-}\lim x=0$, and
$1/4\in\mathrm{L}_x\setminus\Gamma_x(\mathcal{I}_{1/n})$.
Therefore $o(r_n)$ is an optimal rate for both
$\Lambda_x(\mathcal{I}_{1/n})=\mathrm{L}_x$ and
$\Gamma_x(\mathcal{I}_{1/n})=\mathrm{L}_x$, whereas
$O(1/n)$ is optimal within $\mathcal{R}^{\downarrow}$ (for the same properties) 
by Theorem \ref{thm:optimalsummableideal}.
\end{rmk}


\medskip




Now, let $\mathcal{B}$ be the ideal of upper Banach density zero sets, that is, 
$$
\mathcal{B}:=\left\{S\subseteq \omega: \lim_{n\to \infty}\max_{k\in \omega} \frac{|S\cap [k,k+n]|}{n+1}=0\right\},
$$
see e.g. \cite[Definition 2.5]{MR2048515}. 
It is well known that $\mathcal{B}$ is an $F_{\sigma\delta}$ ideal. In addition, it is not a $P$-ideal, see \cite[p. 299]{MR632187}. It follows directly by its definition that $\mathcal{B}\subseteq \mathcal{Z}$ (of course, the inclusion is strict since $\mathcal{Z}$ is a $P$-ideal). Lastly, it is worth noting that $\mathcal{B}$ is not a $P^+$-ideal (and, more generally, all ideals generated by upper densities in the sense of \cite[Definition 1]{MR4054777}, with witnessing decreasing sequence $(2^n\cdot \omega: n \in \omega)$). Hence, thanks to \cite[Proposition 3.2(iii) and Proposition 3.10]{FKL2025}, $\mathcal{B}$ is not a $G_{\delta\sigma}$-ideal (in particular, it is not $F_\sigma$) and there exists a real sequence $x$ such that $\Lambda_x(\mathcal{B})\neq \Gamma_x(\mathcal{B})$. However, as it follows from Theorem \ref{thm:connectedBanach} below, if $x$ satisfies $|\Delta x_n|=o(1)$, then $\Lambda_x(\mathcal{B})=\Gamma_x(\mathcal{B})$. 

\begin{thm}\label{thm:connectedBanach}
Let $\mathcal{I}\neq \mathrm{Fin}$ be an ideal on $\omega$ such that, if $S\subseteq \omega$ contains arbitrarily large finite intervals then $S \in \mathcal{I}^+$ 
\textup{(}which holds, for instance, if $\mathcal{I}\subseteq \mathcal{B}$\textup{)}. 
%

Then $o(1)$ is an optimal rate for both $\Lambda_x(\mathcal{I})=\mathrm{L}_x$ and $\Gamma_x(\mathcal{I})=\mathrm{L}_x$. 
\end{thm}

A nontrivial (non $P^+$) ideal $\mathcal{I}$ contained in $\mathcal{B}$ is the family of sets with upper Buck density zero, that is, more explicitly, the family of all $S\subseteq \omega$ such that $S$ covers $o(n!)$ residues modulo $n!$, see e.g. \cite{MR18196, MR4054777, MR44561}; cf. \cite{MR4439850, MR4746441} for applications. 

We also remark that ideals satisfying the hypothesis of Theorem \ref{thm:connectedBanach} are not necessarily contained in $\mathcal{B}$: this is the case, for instance, of $\mathcal{I}:=\{S\subseteq \omega: |S\cap E|<\infty\}$, where $E:=\{2n: n \in \omega\}$ is the set of even naturals (in fact, $\omega\setminus E\in\mathcal{I}\setminus \mathcal{B}$).


\section{Preliminaries}\label{sec:preliminaries}

Hereafter, we will use the convention $[a,b]:=\omega \cap [a,b]$ for all reals $a,b>0$ (and similarly for other intervals of integers). We start with some preliminary lemmas. 
\begin{lem}\label{lem:technicalbump}
Let $r$ be a nonnegative real sequence such that $\lim_n r_n=0$ and $\sum_n r_n=\infty$.
Fix also $\eta \in \mathbb{R}$, a real $\delta>0$, and integers $i,j\in \omega$ with $i<j$ such that
\begin{equation}\label{eq:conditionsdgkjshd}
\sum_{k=i}^{j-1}r_k>2\delta+2\max\{r_i,r_{i+1},\ldots,r_{j-1}\}.
\end{equation}
Then it is possible to choose reals $x_i,x_{i+1},\ldots,x_j$ such that\textup{:}
\begin{enumerate}[label={\rm (\roman{*})}]
\item $x_i=x_j=\eta$\textup{;}
\item $x_k=\eta+\delta$ for some $k \in (i,j)$\textup{;}
\item $|x_{k+1}-x_k|\le r_k$ for all $k \in [i,j)$\textup{.}
\end{enumerate}
\end{lem}
\begin{proof}
Define $r_\star:=\max\{r_i,r_{i+1},\ldots,r_{j-1}\}$.
Let $t_1\in [i,j)$ be the smallest integer such that
$r_i+r_{i+1}+\cdots+r_{t_1}>\delta$.
Let also $t_2 \in [i,j)$ be the largest integer such that
$r_{j-1}+r_{j-2}+\cdots+r_{t_2}>\delta$.
Note that, thanks to \eqref{eq:conditionsdgkjshd}, $t_1$ and $t_2$
are well defined, and they belong to $[i,j)$.
In addition, we have $t_1<t_2$: in fact, let us suppose for the
sake of contradiction that $t_1\ge t_2$.
Then, by the extremality of such integers, we would obtain
$$
\sum_{k=i}^{j-1}r_k
\le \sum_{k=i}^{t_1}r_k+\sum_{k=t_2}^{j-1}r_k
\le \sum_{k=i}^{t_1-1}r_k+\sum_{k=t_2+1}^{j-1}r_k+2r_\star
\le 2\delta+2r_\star,
$$
where empty sums are understood to be zero.
This contradicts \eqref{eq:conditionsdgkjshd}.
To complete the proof, it is enough to define:
\begin{enumerate}[label={\rm (\roman{*})}]
\item [(a)] $x_k:=\eta+\sum_{h=i}^{k-1}r_h$
for all $k\in[i,t_1]$;
\item [(b)] $x_k:=\eta+\delta$
for all $k\in[t_1+1,t_2]$;
\item [(c)] $x_k:=\eta+\sum_{h=k}^{j-1}r_h$
for all $k\in[t_2+1,j]$.
\end{enumerate}
It follows by construction that the values $x_i,\ldots,x_j$
satisfy the required properties.
\end{proof}

\begin{lem}\label{lem:convergenttauberian}
Let $x$ be a real sequence such that $|\Delta x_n|=o(1)$ and $|\mathrm{L}_x|=1$. Then $x$ is convergent. 
\end{lem}
\begin{proof}
Fix $\eta \in\mathbb{R}$ such that $\mathrm{L}_x=\{\eta\}$. Let us suppose for the sake of contradiction that $x$ is not convergent to $\eta$, hence there exists $\varepsilon>0$ such that $|x_n-\eta|>\varepsilon$ for infinitely many $n \in \omega$. 
Pick a subsequence $(x_{n_k}: k \in \omega)$ such that $|x_{n_{2k}}-\eta|<\varepsilon/4$ and $|x_{n_{2k+1}}-\eta|\ge \varepsilon$ for all $k \in \omega$. 
Since $|\Delta x_n|=o(1)$, this implies that there exists a subsequence $(x_{m_k}: k \in\omega)$ of $x$ taking values in the compact set $K:=\{z \in \mathbb{R}: \varepsilon/3 \le |z-\eta| \le 2\varepsilon/3\}$.  Hence $\mathrm{L}_x \cap K \neq \emptyset$. This contradicts the fact that $\mathrm{L}_x=\{\eta\}$. Therefore $\lim x=\eta$. 
\end{proof}

\medskip

\begin{lem}\label{lem:blocks}
Let $x$ be a real sequence such that $|\Delta x_n|=o(1)$. Then, for each $\eta \in \mathbb{R}$, $\varepsilon>0$, and $m \in \omega$, there exists $n_0 \in \omega$ such that
$$
\forall n\ge n_0, \quad \quad 
|x_n-\eta|<\varepsilon/2 
\quad \implies \quad 
[n,n+m]\subseteq \{k \in \omega: |x_k-\eta|<\varepsilon\}. 
$$
\end{lem}
\begin{proof}
Fix $\eta \in \mathbb{R}$, $\varepsilon>0$, and $m \in \omega$. Choose $n_0 \in \omega$ such that $|\Delta x_n|\le \varepsilon/(2m+1)$ for all $n\ge n_0$. Now, pick $t \in [0,m]$ and $n \ge n_0$ such that $|x_n-\eta|<\varepsilon/2$. Then 
$$
|x_{n+t}-\eta|\le |x_{n+t}-x_n|+|x_n-\eta| \le t\cdot \frac{\varepsilon}{2m+1}+\frac{\varepsilon}{2}<\varepsilon.
$$
Therefore $[n,n+m]\subseteq \{k \in \omega: |x_k-\eta|<\varepsilon\}$.
\end{proof}

\medskip

\begin{lem}\label{lem:partition}
Let $r$ be a positive real sequence such that $\lim_n r_n=0$. Then there exists a partition $(A_k)_{k \in \omega}$ of $\omega$ into finite intervals such that $\lim_k |A_k|=\infty$ and $\lim_k \sum_{n \in A_k}r_n=0$. 
\end{lem}
\begin{proof}
Define $s_n:=\max\{r_k: k\ge n\}$ for each $n \in \omega$, so that the sequence $s=(s_n)_{n \in \omega}$ is decreasing. Now, define recursively the sequence $(n_k)_{k \in \omega}$ in $\omega$ by $n_0:=0$ and 
$$
\forall k \in \omega, \quad 
n_{k+1}:=n_k+\left\lceil \frac{1}{\sqrt{s_{n_k}}}\right\rceil.
$$
To complete the proof, we claim that the intervals $A_k:=[n_k,n_{k+1})$ satisfy the required property. In fact, since $\lim_n s_n=0$, it follows that $\lim_k |A_k|=\lim_k 1/\sqrt{s_{n_k}}=\infty$. Finally, we have 
$$
\forall k \in \omega, \quad 
\sum_{n \in A_k}r_n \le 
|A_k|s_{n_k} \le 
\left(\frac{1}{\sqrt{s_{n_k}}}+1\right)s_{n_k}
=\sqrt{s_{n_k}}+s_{n_k}.
$$
Therefore $\lim_{k} \sum_{n \in A_k}r_n=0$. 
\end{proof}

\medskip

\begin{prop}\label{prop:basicestimateLambdaL}
Let $\mathcal{I}$ be an ideal on $\omega$ and suppose that there exists a partition $(A_k)_{k \in \omega}$ of $\omega$ into finite intervals such that $\bigcup_{k \in J}A_k \in \mathcal{I}^+$ for all infinite $J\subseteq \omega$. Pick also a sequence $r \in \mathcal{R}$ such that 
\begin{equation}\label{eq:conditiononblocks}
\lim_{k\to \infty} \sum_{n \in A_k}r_n=0.
\end{equation}
Then $\Lambda_x(\mathcal{I})=\mathrm{L}_x$ for all real sequences $x$ with $|\Delta x_n|=O(r_n)$. 
\end{prop}
\begin{proof}
Pick a real sequence $x$ and a constant $C>0$ for which $|\Delta x_n|\le C r_n$ for all $n \in \omega$. Since $\Lambda_x(\mathcal{I})\subseteq \mathrm{L}_x$, it will be enough to show the converse inequality. To this aim, suppose that $\mathrm{L}_x$ is nonempty (otherwise the claim is trivial), fix $\eta \in \mathrm{L}_x$, 
and define the infinite sets 
$$
\forall \delta>0, \quad 
S_{\delta}:=\left\{n \in \omega: |x_n-\eta|<\delta \right\}.
$$
Now, pick $n,m,j \in \omega$ and $\delta>0$ such that $n \in S_{\delta/2} \cap A_j$ and $m \in A_j$. It follows that 
\begin{equation}\label{eq:mainestimatelemma}
|x_m-\eta |\le |x_m-x_n|+|x_n-\eta| 
\le \sum_{k \in A_j}|\Delta x_k| + \frac{\delta}{2} \le C\sum_{k \in A_j}r_k + \frac{\delta}{2}.
\end{equation}
Taking into account that $\lim_j\sum_{k \in A_j}r_k=0$, it follows that $|x_m-\eta| < \delta$ for every sufficiently large $j$ with $S_{\delta/2}\cap A_j\neq\emptyset$ and every $m \in A_j$. Hence $S_\delta$ contains infinitely many intervals $A_j$, i.e., 
\begin{equation}\label{eq:claiminclusiobfkjhksfjhg}
\exists j^\star_\delta \in\omega, \quad 
\bigcup_{j \ge j^\star_\delta, \,\, S_{\delta/2} \cap A_j \neq \emptyset} A_j \subseteq S_\delta. 
\end{equation}

This implies that there exists a strictly increasing sequence $(j_m: m \in \omega)$ in $\omega$ such that $j_m \ge j_{2^{-m}}^\star$ and $S_{2^{-m-1}} \cap A_{j_m}\neq \emptyset$. To complete the proof, it is enough to show that 
the subsequence $(x_n: n \in A)$ is convergent to $\eta$, where $A:=\bigcup_m A_{j_m}$ (of course, $A \in \mathcal{I}^+$ since it contains infinitely many $A_j$). In fact, for each $k \in \omega$, we obtain by \eqref{eq:claiminclusiobfkjhksfjhg}
$$
\{n \in A: |x_n-\eta|<2^{-k}\}=S_{2^{-k}}\cap A\supseteq \bigcup_{m\ge k}A_{j_m}.
$$
By the arbitrariness of $k$, we conclude that $\eta \in \Lambda_x(\mathcal{I})$. Therefore $\Lambda_x(\mathcal{I})=\mathrm{L}_x$.
\end{proof}
\begin{rmk}
Let $(A_k)_{k \in \omega}$ be a sequence as in the statement of Proposition \ref{prop:basicestimateLambdaL}. 
It follows by Theorem \ref{equivalencefridyclaim} and Proposition \ref{prop:basicestimateLambdaL} that, if $r \in \mathcal{R}$ satisfies condition \eqref{eq:conditiononblocks}, then $O(r_n)$ is a Tauberian condition for $\Gamma_x(\mathcal{I})=\mathrm{L}_x$. 
%
%

However, the converse fails (namely, condition \eqref{eq:conditiononblocks} is sufficient but not necessary for the latter conclusion): in fact, suppose that $\mathcal{I}=\mathcal{Z}$, and define $A_0:=\{0,1\}$ and $A_k:=[2^k,2^{k+1})$ for all $k\ge 1$. It follows by Fridy's Theorem \ref{thm:fridythm3} and Theorem \ref{equivalencefridyclaim} that $\Gamma_x(\mathcal{Z})=\mathrm{L}_x$ for all real sequences $x$ with $|\Delta x_n|=O(1/n)$. On the other hand, it is easily seen that $\lim_k \sum_{n \in A_k}1/n=\log(2)$, hence condition \eqref{eq:conditiononblocks} does not hold. 
\end{rmk}


\medskip

\begin{prop}\label{prop:basic2estimateLambdaLwitho}
Let $\mathcal{I}$ and $(A_k)_{k \in \omega}$ as in the statement of Proposition \ref{prop:basicestimateLambdaL}. Pick also a sequence $r \in \mathcal{R}$ such that 
$$
\sup_{k\in \omega} \sum_{n \in A_k}r_n<\infty. 
$$
Then $\Lambda_x(\mathcal{I})=\mathrm{L}_x$ for all real sequences $x$ with $|\Delta x_n|=o(r_n)$. 
\end{prop}
\begin{proof}
It goes verbatim as in the proof of Proposition \ref{prop:basicestimateLambdaL}, replacing the last estimate in \eqref{eq:mainestimatelemma} with 
$$
\sum_{k \in A_j}|\Delta x_k|+\frac{\delta}{2} \le \left(\max_{k \in A_j}\frac{|\Delta x_k|}{r_k}\right) \cdot \sup_{t\in \omega}\sum_{n \in A_t}r_n +\frac{\delta}{2}=o(1)+\frac{\delta}{2}.
$$
We omit further details. 
\end{proof}

\section{Proofs of Main Results}\label{sec:proofs}

Now, we are ready to prove our main results. 

\begin{proof}[Proof of Theorem \ref{equivalencefridyclaim}]
\ref{item:A1} $\implies$ \ref{item:A2}. This is clear since $\Lambda_x(\mathcal{I})\subseteq \Gamma_x(\mathcal{I})\subseteq \mathrm{L}_x$. 

\medskip

\ref{item:A2} $\implies$ \ref{item:A3}. Pick a real sequence $x=(x_n)_{n \in \omega}$ and a constant $C>0$ such that $|\Delta x_n|\le Cr_n$ for all $n \in \omega$.  Suppose also that $\mathcal{I}\text{-}\lim x=\eta$ for some $\eta \in \mathbb{R}$, so that $\Gamma_x(\mathcal{I})=\{\eta\}$. Thanks to item \ref{item:A2}, we have $\mathrm{L}_x=\{\eta\}$. We need to show that $\lim x=\eta$. 

To this aim, let us suppose for the sake of contradiction that there exists $\varepsilon>0$ and infinitely many $n \in \omega$ such that $|x_n-\eta|>\varepsilon$. By $\mathcal{I}$-convergence, there are also infinitely many $n \in \omega$ such that $|x_n-\eta|<\varepsilon/4$. Hence there exists a subsequence $(x_{n_k}: k \in \omega)$ such that $|x_{n_{2k}}-\eta|<\varepsilon/4$ and $|x_{n_{2k+1}}-\eta|\ge \varepsilon$ for all $k \in \omega$.  
The conclusion follows 
with the same argument as that given at the end
of the proof of Lemma \ref{lem:convergenttauberian}.

\medskip

\ref{item:A3} $\implies$ \ref{item:A2}. Pick a real sequence $x=(x_n)_{n \in \omega}$ and a constant $C>0$ such that $|\Delta x_n|\le Cr_n$ for all $n \in \omega$. It will be enough to show that $\mathrm{L}_x\subseteq \Gamma_x(\mathcal{I})$. This is clear if $\mathrm{L}_x=\emptyset$. Otherwise, fix $\eta \in \mathrm{L}_x$. Suppose for the sake of contradiction that there exists $\varepsilon>0$ such that $S:=\{n \in \omega: |x_n-\eta|< \varepsilon\} \in \mathcal{I}$. In particular, $S$ is not cofinite. Since $\eta$ is an accumulation point, the set $\{n \in \omega: |x_n-\eta|<\varepsilon/2\}$ is infinite. 
This implies there exists a subsequence $(x_{n_k}: k \in \omega)$ such that $|x_{n_{2k}}-\eta|< \varepsilon/2$ and $|x_{n_{2k+1}}-\eta|\ge \varepsilon$ for all $k \in \omega$. 
Define $\varepsilon^\prime:=\varepsilon/8C$ for notational convenience. 
Also, without loss of generality, it can be assumed that $n_0>0$, that $r_n<\varepsilon^\prime/2$ for all $n\ge n_0$, and that the index $n_{2k+1}$ is the minimal integer (greater than $n_{2k}$) for which $|x_{n_{2k+1}}-\eta|\ge \varepsilon$. Define the interval $J_k:=[n_{2k}, n_{2k+1})$ for each $k \in \omega$, so that $J:=\bigcup_k J_k \subseteq S\in \mathcal{I}$. In addition, for each $k \in \omega$ we have 
$$
\frac{\varepsilon}{2} \le |x_{n_{2k+1}}-x_{n_{2k}}| \le \sum_{n \in J_k}|\Delta x_n| \le C\sum_{n \in J_k}r_n.
$$
This implies that $\sum_{n \in J_k}r_n \ge 4\varepsilon^\prime$ for all $k \in \omega$. 

Now, let us construct a new sequence $y$ as follows: $y_n:=0$ for all $n \notin J$; otherwise, pick $k \in \omega$, and observe that 
$$
\sum_{t=n_{2k}}^{n_{2k+1}-2}r_t 
> -\frac{\varepsilon^\prime}{2}+\sum_{t \in J_k}r_t 
> 3\varepsilon^\prime>2\varepsilon^\prime+2\max\{r_t: t \in J_k\}.
$$
It follows by Lemma \ref{lem:technicalbump} that, on each interval $J_k=[n_{2k},n_{2k+1})$, 
it is possible to select values $(y_n: n \in J_k)$ such that: $y_{n_{2k}}=y_{n_{2k+1}-1}=0$, $|\Delta y_n|\le r_n$ for all $n \in J_k$, and $y_n=\varepsilon^\prime$ for some $n \in J_k$. 
In addition, $\{n \in \omega: y_n \neq 0\}\subseteq J\subseteq S$, so that $\mathcal{I}\text{-}\lim y=0$. It follows by item \ref{item:A3} that $\lim y=0$. This contradicts the fact that $y_n=\varepsilon^\prime$ for infinitely many $n$. Therefore $\eta$ has to be an $\mathcal{I}$-cluster point of $x$. 

\medskip

To prove the second part, suppose that $\mathcal{I}$ is a $P^+$-ideal. We need to only show that the converse implication \ref{item:A2} $\implies$ \ref{item:A1}. Pick a real sequence $x$ such that $|\Delta x_n|=O(r_n)$. Then it follows by \cite[Theorem 3.4]{MR4393937} that $\Lambda_x(\mathcal{I})=\Gamma_x(\mathcal{I})$ for all real sequences $x$; cf. also \cite[Proposition 3.10]{FKL2025}, and 
\cite[Theorem 2.3]{PaoloMarek17} or 
\cite[Proposition 5.4]{MR4967065}
for the special case of $F_\sigma$ ideals. This concludes the proof.
\end{proof}

\medskip

\begin{proof}
[Proof of Theorem \ref{thm:charactbaire}]
\ref{item:BP1} $\implies$ \ref{item:BP2A}. 
Suppose that $\mathcal{I}$ has the Baire property. Hence it is possible to fix a partition $(A_j: j \in \omega)$ of $\omega$ such that each $A_j$ is a finite interval of nonnegative integers and $S\notin \mathcal{I}$ whenever $A_j \subseteq S$ for infinitely many $j \in \omega$; see \cite[Theorem 2.1]{MR579439} and also \cite[Theorem 4.1.2]{MR1350295}. Fix an arbitrary sequence of positive reals $q=(q_j: j \in \omega)$ such that $\lim_j q_j=0$, and $\sum_j q_j=\infty$ (i.e., $q \in \mathcal{R}$), and define 
\begin{equation}\label{eq:kjdshgsd}
r_n:=\frac{q_j}{|A_j|}
\end{equation}
for all $n,j \in \omega$ such that $n \in A_j$ (for instance, one can pick $q_j:=1/(j+1)$ for all $j \in \omega$). 
Observe that $\lim_j \sum_{n \in A_j}r_n=0$. It follows by Proposition \ref{prop:basicestimateLambdaL} that $\Lambda_x(\mathcal{I})=\mathrm{L}_x$ for all real sequences $x$ with $|\Delta x_n|=O(r_n)$. 

\medskip

\ref{item:BP2A} $\implies$ \ref{item:BP2}. This is clear. 

\medskip


\ref{item:BP2} $\implies$ \ref{item:BP3}. This follows by 
Theorem \ref{equivalencefridyclaim}. 

\medskip

\ref{item:BP3} $\implies$ \ref{item:BP1}. 
Fix a positive real sequence $(r_n)_{n \in \omega}$ satisfying the property in item \ref{item:BP3}. Since $\lim_nr_n=0$ and $\sum_n r_n=\infty$, it is possible to pick a partition $(A_j)_{j \in \omega}$ of $\omega$ into finite intervals of nonnegative integers such that $\sum_{n \in A_j}r_n >4$ for all $j \in \omega$ and $r_n < 1/2$ for all $n\ge \min A_1$. Now, suppose for the sake of contradiction that $\mathcal{I}$ does not have the Baire property. Then there exists an infinite set $J\subseteq \omega\setminus \{0\}$ such that $Q:=\bigcup_{j \in J}A_j \in\mathcal{I}$. Notice that $Q$ cannot be cofinite. At this point, let us construct a real sequence $(x_n)_{n \in \omega}$ as  follows: 
If $j\notin J$ and $n \in A_j$ then $x_n:=0$.

Otherwise, pick $j \in J$, set $\delta:=1$ and observe that
$$
\sum_{n=\min A_j}^{\max A_j-1}r_n
=\sum_{n \in A_j}r_n-r_{\max A_j}
>4-\frac12
>2\delta+2\max\{r_n: n \in A_j\}.
$$
It follows by Lemma \ref{lem:technicalbump}, with $\eta:=0$, that,
on each interval $A_j$, it is possible to select values
$(x_n: n \in A_j)$ such that:
$x_{\min A_j}=x_{\max A_j}=0$,
$|\Delta x_n|\le r_n$ for all $n \in A_j\setminus\{\max A_j\}$,
and $x_n=1$ for some $n \in A_j$.
Since the increments across block boundaries are zero, the bound
$|\Delta x_n|\le r_n$ also holds for $n=\max A_j$.
It follows by construction that $\{0,1\}\subseteq \mathrm{L}_x$, hence $x$ is not convergent. 
On the other hand, it is clear that $\{n \in \omega: x_n \neq 0\}\subseteq Q \in \mathcal{I}$. It follows that $\mathcal{I}\text{-}\lim x=0$, which 
completes the proof. 
\end{proof}

\medskip

\begin{proof}
[Proof of Theorem \ref{thm:conjectureweak}] 
Since $\mathcal{Z}$ is a Baire ideal with witnessing partition $(A_j:j \in \omega)$ defined by $A_0:=\{0,1\}$ and $A_j:=[2^{j},2^{j+1})$,  
it follows by Theorem \ref{thm:charactbaire} and Remark \ref{rmk:choicern} that items \ref{item:BP2A}, \ref{item:BP2}, and \ref{item:BP3} hold with 
$$
r_n=o\left(\frac{1}{2^{1+\lfloor \log_2(n)\rfloor}-2^{\lfloor \log_2(n)\rfloor}}\right)=o\left(\frac{1}{n}\right). 
$$
This completes the proof. 
\end{proof}

\medskip

\begin{proof}
[Proof of Proposition \ref{prop:Baire}]
Thanks to Lemma \ref{lem:partition}, it is possible to pick a partition $(A_k)_{k \in \omega}$ of $\omega$ into finite intervals such that $\lim_k |A_k|=\infty$ for all $k \in \omega$ and $\lim_k\sum_{n \in A_k}r_n=0$. At this point, it is enough to show that the family 
$$
\mathcal{I}:=\left\{S\subseteq \omega: \sup_{k\in \omega} |S\cap A_k|<\infty\right\}
$$
satisfies our claim. First of all, it is clear that $\mathcal{I}$ is an ideal and 
that it 
contains the infinite set $\bigcup_k \{\min A_k\}$. In addition, it can be written as $\bigcup_{m \in \omega}\mathcal{I}_m$, where each $\mathcal{I}_m:=\bigcap_{k \in \omega}\{S\subseteq \omega: |S\cap A_k|\le m\}$ is closed. Hence $\mathcal{I}$ is an $F_\sigma$ ideal on $\omega$. 
It is also easy to see that it is not a $P$-ideal by picking any witnessing sequence $(S_n)_{n \in \omega}$ such that $|S_n\cap A_k|\le \min\{n,|A_k|\}$ for all $n,k \in \omega$. 

Observe that, if $S\subseteq \omega$ contains infinitely many $A_k$, then $S \in \mathcal{I}^+$. Hence $(A_k)_{k \in \omega}$ is a witnessing sequence for the Baire property of $\mathcal{I}$. We conclude by Proposition  \ref{prop:basicestimateLambdaL} that $\Lambda_x(\mathcal{I})=\mathrm{L}_x$ for all real sequences $x$ with $|\Delta x_n|=O(r_n)$. 
\end{proof}

\medskip

\begin{proof}
[Proof of Theorem \ref{thm:mainanalyticPideal}] 
 \ref{item:Analytic1} $\Longleftrightarrow$  \ref{item:Analytic1bis}. It follows by the equivalence \ref{item:A2} $\Longleftrightarrow$  \ref{item:A3} in Theorem \ref{equivalencefridyclaim}.

 \medskip
 
 \ref{item:Analytic1bis} $\implies$  \ref{item:Analytic2}. Let us suppose that $\varphi: \mathcal{P}(\omega)\to [0,\infty]$ is a lscsm such that $\mathcal{I}=\mathrm{Exh}(\varphi)$ and $\varphi$ is not smooth. We claim that there exists a non-convergent real sequence $(x_n)_{n \in \omega}$ such that $|\Delta x_n|\le 1/n$ for all $n\ge 1$ and $\mathcal{I}\text{-}\lim_n x_n=0$.

 
 Since $\varphi$ is not smooth, it is possible to fix a real $\kappa>1$ such that $\liminf_n \varphi([n,\kappa n])=0$. Hence, there exists an increasing sequence of nonnegative integers $(n_t)_{t \in \omega}$ such that $\lim_t \varphi([n_t,\kappa n_t])=0$. Without loss of generality, we can suppose also that $n_0>0$, $n_{t+1}>1+\kappa n_t$, $\kappa^{1/3}n_t>n_t+2$, and $\varphi(W_t)<2^{-t}$ for all $t \in \omega$, where $W_t:=[n_t, \kappa n_t]$. 
 
 At this point, define
$$
\eta_0:=\frac{\log(\kappa)}{6}.
$$
Since
$$
\lim_{t\to\infty}
\sum_{n=n_t}^{\lfloor\kappa n_t\rfloor-1}\frac{1}{n}
=\log(\kappa)>2\eta_0,
$$
after passing to a tail and relabelling, we may suppose that
$$
\sum_{n=n_t}^{\lfloor\kappa n_t\rfloor-1}\frac{1}{n}
>2\eta_0+\frac{2}{n_t}
$$
for all $t\in\omega$.
It follows by Lemma \ref{lem:technicalbump}, with $\eta:=0$,
$\delta:=\eta_0$, and increment bounds $1/n$, that on each
interval $W_t$ it is possible to select values $(x_n:n\in W_t)$
such that:
\begin{enumerate}[label={\rm (\arabic{*})}]
\item $x_{n_t}=x_{\lfloor\kappa n_t\rfloor}=0$;
\item $x_n=\eta_0$ for some $n\in W_t$;
\item $|\Delta x_n|\le 1/n$ for all
$n\in[n_t,\lfloor\kappa n_t\rfloor)$.
\end{enumerate}
Set $x_n:=0$ for all $n\in\omega\setminus\bigcup_t W_t$.
Since the intervals are separated and their endpoint values
are zero, we obtain $|\Delta x_n|\le 1/n$ for all $n\ge1$.

By construction, $\eta_0\in\mathrm{L}_x$ and $x_{n_t}=0$
for all $t\in\omega$. Since $\eta_0>0$, the sequence $x$
is not convergent. On the other hand, we have 
$$
\{n \in \omega: x_n\neq 0\}\subseteq W:=\bigcup_{t \in \omega}W_t
$$
and, at the same time, $W$ belongs to $\mathcal{I}$ since 
\begin{displaymath}
\begin{split}
\left\|W\right\|_\varphi&=\lim_{t\to \infty}\varphi(W\setminus [0,n_t)) 
=\lim_{t\to \infty} \lim_{y\to \infty}\varphi(W\cap [n_t,y])\\
&\le \lim_{t\to \infty} \lim_{i\to \infty}\varphi(W_t \cup \cdots \cup W_{t+i}) 
\le \lim_{t\to \infty} \sum_{j\ge t} \varphi(W_j)
\le \lim_{t\to \infty} \sum_{j\ge t} 2^{-j} =0.
\end{split}
\end{displaymath}
Therefore $\{n \in \omega: x_n\neq 0\} \in \mathcal{I}$, so that $\mathcal{I}\text{-}\lim_n x_n=0$. 
 
\medskip

 \ref{item:Analytic2} $\implies$  \ref{item:Analytic3}. This is clear. 
 
 \medskip

 \ref{item:Analytic3} $\implies$  \ref{item:Analytic1}. 
Let $\varphi: \mathcal{P}(\omega)\to [0,\infty]$ be a smooth lscsm such that $\mathcal{I}=\mathrm{Exh}(\varphi)$. 
Let $(x_n)_{n \in \omega}$ be a real sequence such that $|\Delta x_n|=O(1/n)$ as $n\to \infty$. We need to show that $\Gamma_x(\mathcal{I})=\mathrm{L}_x$. 
To this aim, since $\Gamma_x(\mathcal{I})\subseteq \mathrm{L}_x$, we need to show that the converse inequality holds as well. Suppose that $\mathrm{L}_x\neq \emptyset$ (otherwise there is nothing to show), and fix $\eta \in \mathrm{L}_x$. Fix also a constant $C>0$ such that $|\Delta x_n| \le C/n$ for all $n\ge 1$. 

Pick a subsequence $(x_{n_t}: t \in \omega)$ which is convergent to $\eta$. Fix $\varepsilon>0$ and pick a real 
\begin{equation}\label{eq:condition1kappa}
\kappa  \in \left(1,\,e^{\frac{\varepsilon}{4C}}\right).
\end{equation} 
At this point, suppose without loss of generality that the subsequence is sufficiently sparse, namely, $\kappa n_0>n_0+1$ and that $n_{t+1}>1+\kappa n_t$ for all $t \in \omega$. This implies that the sequence $(W_t)_{t \in \omega}$ defined by $W_t:=\omega \cap [n_t, \kappa n_t]$ consists of  nonempty sets of nonnegative integers with $\max W_t<\min W_{t+1}$ for all $t \in \omega$. Fix a sufficiently large integer $t_0 \in \omega$ such that $|x_{n_t}-\eta|<\varepsilon/2$ for all $t\ge t_0$ and 
\begin{equation}\label{eq:condition2kappa}
4C+\varepsilon<\varepsilon n_{t_0}.
\end{equation}
It follows that, if $n \in W_{t}$ for some $t\ge t_0$, then 
\begin{displaymath}
\begin{split}
|x_n-x_{n_t}| &\le \sum_{i=0}^{n-n_t-1}|\Delta x_{n_t+i}| 
\le C\sum_{j=n_t}^{n-1} \frac{1}{j}
\le C\int_{n_t}^n \frac{1}{x-1}\mathrm{d}x
\\&
\le C \log\left(\frac{\kappa n_t}{n_t-1}\right)
\le C\log(\kappa)+C\log\left(1+\frac{1}{n_{t_0}-1}\right)
\\&
<\frac{\varepsilon}{4}+\frac{C}{n_{t_0}-1}
<\frac{\varepsilon}{4}+\frac{\varepsilon}{4}=\frac{\varepsilon}{2}.
\end{split}
\end{displaymath}
In the last upper estimates above, we used \eqref{eq:condition1kappa} and \eqref{eq:condition2kappa}. 
Hence, for all integers $n \in W_{t}$ with $t\ge t_0$, we get 
\begin{displaymath}
|x_n-\eta| \le |x_n-x_{n_t}|+|x_{n_t}-\eta| 
< \frac{\varepsilon}{2}+\frac{\varepsilon}{2}=\varepsilon.
\end{displaymath}

It follows that the set $\{n \in \omega: |x_n-\eta|<\varepsilon\}$ contains $W_{t}$ for all $t\ge t_0$. Taking into account that $\varphi$ is a smooth lscsm, we conclude that 
$$
\left\|\{n \in \omega: |x_n-\eta|<\varepsilon\}\right\|_\varphi
\ge \left\|\bigcup_{t\ge t_0}\,W_{t}\right\|_\varphi
\ge \limsup_{t\to \infty}\varphi(W_{t})>0. 
$$
Therefore $\{n \in \omega: |x_n-\eta|<\varepsilon\}\notin \mathcal{I}$. By the arbitrariness of $\varepsilon$, we get $\eta \in \Gamma_x(\mathcal{I})$. 
\end{proof}

\medskip

\begin{proof}
[Proof of Corollary \ref{cor:originalconjecture}] 
It is sufficient to prove the claim for the ideal $\mathcal{Z}$. 
Define the lscsm $\varphi: \mathcal{P}(\omega) \to [0,\infty]$ by 
$$
\forall S\subseteq \omega, \quad 
\varphi(S):=\sup_{n \in \omega}\frac{|S\cap [0,n]|}{n+1}.
$$
Then it is clear that $\|S\|_\varphi$ is the upper asymptotic density of $S\subseteq \omega$, and $\mathcal{Z}=\mathrm{Exh}(\varphi)$. At this point, thanks to Theorem \ref{thm:mainanalyticPideal}, it is enough to show that   $\varphi$ is smooth. In fact, it is easy to see that 
$$
\forall \kappa>1, \quad 
\liminf_{n\to \infty} \varphi([n,\kappa n])
\ge \liminf_{n\to \infty} \frac{|[n, \kappa n]|}{1+\lfloor \kappa n\rfloor} 
= \frac{\kappa-1}{\kappa}>0.
$$
This completes the proof. 
\end{proof}

\medskip

\begin{proof}
[Proof of Corollary \ref{cor:recovercorollary}] 
As in the previous proof, it is sufficient to prove the claim for the ideal $\mathcal{Z}$. Let $x$ be a (not necessarily bounded) real sequence, and fix $\eta \in \mathbb{R}$. For each $\varepsilon>0$, define the bounded real sequence $y^\varepsilon$ by 
$$
\forall n \in \omega, \quad 
y^\varepsilon_n:=\lambda_n\left(\left\{k \in \omega: |x_k-\eta|>\varepsilon \right\}\right).
$$
Observe that $|y^\varepsilon_{n+1}-y^\varepsilon_n|=O(1/n)$ as $n\to \infty$. Now, the sequence $x$ is $\mathcal{Z}$-statistically convergent to $\eta$ if and only if $\mathcal{Z}\text{-}\lim y^\varepsilon=0$ for all $\varepsilon>0$, which is equivalent by Corollary  \ref{cor:originalconjecture} to $\lim y^\varepsilon=0$ for all $\varepsilon>0$. The latter is simply a rewriting of $\mathcal{Z}\text{-}\lim x =\eta$. 
\end{proof}

\medskip

\begin{proof}
[Proof of Proposition \ref{prop:summableborderline}] 
Set $P:=\{2^k: k \in \omega\}$ and let $a$ be the nonnegative real sequence defined by $a_n=1$ if $n \in P$ and $a_n:=0$ otherwise. We claim that the summable ideal $\mathcal{I}$ defined as in \eqref{eq:summable} has the required properties. For, observe that 
$$
\mathcal{I}=\{S\subseteq \omega: |S\cap P|<\infty\}.
$$
In other words, $\mathcal{I}$ is a copy on $\omega$ of the Fubini sum $\mathrm{Fin}\oplus \mathcal{P}(\omega)$, see e.g. \cite[Chapter 1]{MR1711328}. 

\smallskip

Define the interval of integers $A_j:=[2^j,2^{j+1})$ for each $j \in \omega$. Pick an infinite set $J\subseteq \omega$ and suppose that $S\subseteq \omega$ contains $\bigcup_{j \in J}A_j$. Then $S$ contains infinitely many elements of $P$, so that $S \notin \mathcal{I}$. Then item \ref{item:example1fixed} follows by Remark \ref{rmk:o1overnBaire}. 

\smallskip

Now, let $\varphi: \mathcal{P}(\omega)\to [0,\infty]$ be the lscsm defined by $\varphi(S):= \sum_{n \in S}a_n$ for all $S\subseteq \omega$. 
It is easy to see that $\mathcal{I}=\mathrm{Exh}(\varphi)$. In addition, $\varphi$ is not smooth: in fact, setting $\kappa=3/2$ in \eqref{eq:definitionsmooth}, it follows that, by the monotonicity of $\varphi$, 
$$
\varphi([2^j+1,\, \kappa (2^j+1)])\le 
\varphi([2^j+1,2^{j+1}-1])=0
$$
for all sufficiently large integers $j$. This implies that $\liminf_n \varphi([n,\kappa n])=0$. Thanks to Theorem \ref{thm:mainanalyticPideal} and the fact that there exists a non-smooth lscsm $\varphi$ for which $\mathcal{I}=\mathrm{Exh}(\varphi)$, items \ref{item:example2fixed} and \ref{item:example3fixed} follow. 
\end{proof}

\medskip

\begin{proof}
[Proof of Proposition \ref{prop:ratesummableo}] 
Up to rescaling the sequence $r$ by a positive constant, suppose without loss of generality that $r_n \le 1$ for all $n \in \omega$. Define recursively a strictly increasing sequence $(n_k)_{k \in \omega}$ in $\omega$ such that $n_0:=0$ and, for all $k \in \omega$: 
\begin{enumerate}[label={\rm (\roman{*})}]
\item $n_{2k+1}$ is the smallest integer $t>n_{2k}$ such that $\sum_{i=n_{2k}}^{t-1} r_i \ge 6$; 
\item $n_{2k+2}$ is the smallest integer $t>n_{2k+1}$ such that $\sum_{i=n_{2k+1}}^{t-1} r_i \ge 1$. 
\end{enumerate}
Set $A_k:=[n_k,n_{k+1})$ for all $k \in \omega$, and define the summable ideal
$$
\mathcal{I}:=\left\{S\subseteq \omega: \sum_{k \in \omega}\sum_{n \in S\cap A_{2k+1}}\frac{1}{|A_{2k+1}|}<\infty\right\}.
$$
Of course, the sequence $(B_k)_{k \in \omega}$, where $B_k:=A_{2k}\cup A_{2k+1}$, witnesses the Baire property of $\mathcal{I}$. To complete the proof, we claim that $\mathcal{I}$ satisfies the claimed property. 

To this aim, pick a real sequence $x$ such that $|\Delta x_n|=o(r_n)$ and observe that 
$$
\sup_{k\in \omega} \sum_{n \in B_k}r_n \le \sup_{k \in \omega}\left(7+2\max\{r_n: n \in B_k\}\right) \le 9.
$$
It follows by Proposition \ref{prop:basic2estimateLambdaLwitho} that $\Lambda_x(\mathcal{I})=\mathrm{L}_x$, hence also $\Gamma_x(\mathcal{I})=\mathrm{L}_x$. 

Vice versa, observe that 
$$
\sum_{n \in [n_{2k},n_{2k+1}-2]}r_n \ge -1+\sum_{n \in A_{2k}}r_n \ge -1+6>2\cdot 1+2\max\{r_n: n \in \omega\}.
$$
It follows by Lemma \ref{lem:technicalbump} that there exists a real sequence $x$ supported on $A:=\bigcup_k A_{2k}$ which satisfies the following properties: 
\begin{enumerate}[label={\rm (\roman{*})}] 
\item $x_n=0$ for all $n \notin A$, 
\item $x_{\min A_{2k}}=x_{\max A_{2k}}=0$ for all $k \in \omega$, 
\item for each $k \in \omega$ there exists $n \in A_{2k}$ such that $x_n=1$, and 
\item $|\Delta x_n|\le r_n$ for all $n \in \omega$. 
\end{enumerate}
This implies that $\{0,1\}\subseteq \mathrm{L}_x$. On the other hand, $\{n \in \omega: x_n\neq 0\}\subseteq A\in \mathcal{I}$. Taking into account that $\mathcal{I}$ is a $P^+$ ideal (so that $\Lambda_x(\mathcal{I})=\Gamma_x(\mathcal{I})$), this implies that $1 \in \mathrm{L}_x\setminus \Gamma_x(\mathcal{I})$.
\end{proof}

\medskip

\begin{proof}
[Proof of Theorem \ref{thm:optimalsummableideal}]
First, since $\mathcal{I}\subseteq \mathcal{Z}$, it follows by Corollary \ref{cor:originalconjecture} that $O(1/n)$ is a Tauberian condition for $\Gamma_x(\mathcal{I})=\mathrm{L}_x$. Vice versa, fix a sequence $r^\prime \in \mathcal{R}^\downarrow$ and suppose that $r^\prime\neq O(1/n)$, that is, $\sup_nnr^\prime_n=\infty$. Then it is possible to pick a strictly increasing sequence $(n_k)_{k \in \omega}$ of positive integers such that $r^\prime_{n_0}<1$ and, in addition, 
\begin{equation}\label{eq:disjointsummable}
\forall k \in \omega, \quad 2n_k r^\prime_{2n_k}>2^{k+4} 
\quad \text{ and }\quad 
n_{k+1}>4 n_k.
\end{equation}
Define the interval $I_k:=[n_k, n_k+8/r^\prime_{2n_k})$ for each $k \in \omega$, and observe that $I_k\subseteq [n_k,2n_k)$ by \eqref{eq:disjointsummable}. Now, we have by construction
$$
\forall k \in \omega, \quad 
\sum_{n \in I_k}r^\prime_n \ge |I_k|r^\prime_{\max I_k}\ge \left(\frac{8}{r^\prime_{2n_k}}-1\right)\cdot r^\prime_{2n_k}>7>3+2 r^\prime_{n_0}.
$$
Hence $\sum_{n=\min I_k}^{\max I_k-1}r^\prime_n>2\cdot 1+2 r^\prime_{n_0}$. 
It follows by Lemma \ref{lem:technicalbump} that there exists a real sequence $x$ supported on $I:=\bigcup_k I_k$ which satisfies the following properties: 
\begin{enumerate}[label={\rm (\roman{*})}] \item 
$x_n=0$ for all $n \notin I$, 
\item $x_{\min I_k}=x_{\max I_k}=0$ for all $k \in \omega$, 
\item for each $k \in \omega$ there exists $n \in I_k$ such that $x_n=1$, and 
\item $|\Delta x_n|\le r^\prime_n$ for all $n \in \omega$. 
\end{enumerate}
This implies that $\{0,1\}\subseteq \mathrm{L}_x$. To complete the proof, it will be enough to show that $1$ is not an $\mathcal{I}$-cluster point of $x$. Taking into account \eqref{eq:disjointsummable}, it follows that 
$$
\sum_{n \in I_k}\frac{1}{n}
\le 2\log\left(\frac{n_k+\frac{8}{r^\prime_{2n_k}}}{n_k}\right)
\le 2\cdot \frac{8}{n_kr^\prime_{2n_k}}\le \frac{32}{2^{k+4}}
$$
for all sufficiently large $k$. This implies that $\sum_{n \in I}1/n \le O(1)+\sum_k 2^{1-k}<\infty$, hence $I \in \mathcal{I}_{1/n}\subseteq \mathcal{I}$. Therefore $\mathcal{I}\text{-}\lim x=0$ and, in particular, $1\notin \Gamma_x(\mathcal{I})$. 
\end{proof}

\medskip

\begin{proof}
[Proof of Theorem \ref{thm:connectedBanach}] 
%
Let $x$ be a real sequence such that $|\Delta x_n|=o(1)$. 
First, let us show that $\mathrm{L}_x\subseteq \Lambda_x(\mathcal{I})$. This is clear if $\mathrm{L}_x=\emptyset$. Otherwise, fix an accumulation point $\eta \in \mathrm{L}_x$. Thanks to Lemma \ref{lem:blocks}, it is possible to define recursively a sequence $(n_k)_{k \in \omega}$ in $\omega$ such that: \begin{enumerate}[label={\rm (\roman{*})}]
\item $n_{k+1}>n_k+k$ for all $k \in \omega$;
\item $|x_n-\eta|<2^{-k}$ for all $k \in \omega$ and $n \in [n_k,n_k+k]$. 
\end{enumerate}
It follows by construction that $S:=\bigcup_k [n_k,n_k+k]$ is an $\mathcal{I}$-positive set, and the subsequence $(x_n: n \in S)$ is convergent to $\eta$. Therefore $\Lambda_x(\mathcal{I})=\Gamma_x(\mathcal{I})=\mathrm{L}_x$. The conclusion on the optimal rate $o(1)$ follows by Remark \ref{rmk:IFin}.
\end{proof}


\bibliographystyle{amsplain}
\bibliography{idealem}

\end{document}